\documentclass{amsproc}
\usepackage{amssymb}
\usepackage{amsmath}
\usepackage{eurosym}
\usepackage{amsfonts}

            \advance \topmargin by -\headheight \advance
            \topmargin by -\headsep
            \evensidemargin
            \oddsidemargin
\newtheorem{theorem}{Theorem}[section]
\newtheorem{lemma}[theorem]{Lemma}

\newtheorem{remark}[theorem]{Remark}
\numberwithin{equation}{section}
\input{tcilatex}

\begin{document}
\title[Tribonacci and Tribonacci-Lucas Sedenions]{Tribonacci and
Tribonacci-Lucas Sedenions}
\thanks{}
\author[Y\"{u}ksel ~Soykan]{Y\"{U}KSEL SOYKAN}
\maketitle

\begin{center}
\textsl{Zonguldak B\"{u}lent Ecevit University, Department of Mathematics, }

\textsl{Art and Science Faculty, 67100, Zonguldak, Turkey }

\textsl{e-mail: \ yuksel\_soykan@hotmail.com}
\end{center}

\textbf{Abstract.} The sedenions form a 16-dimensional Cayley-Dickson
algebra. In this paper, we introduce the Tribonacci and Tribonacci-Lucas
sedenions. Furthermore, we present some properties of these sedenions and
derive relationships between them.



\section{Introduction}

Tribonacci sequence $\{T_{n}\}_{n\geq 0}$ and Tribonacci-Lucas sequence $%
\{K_{n}\}_{n\geq 0}$ are defined by the third-order recurrence relations%
\begin{equation}
T_{n}=T_{n-1}+T_{n-2}+T_{n-3},\text{ \ \ \ \ }T_{0}=0,T_{1}=1,T_{2}=1,
\label{equati:fvcvxghsbnz}
\end{equation}%
and 
\begin{equation}
K_{n}=K_{n-1}+K_{n-2}+K_{n-3},\text{ \ \ \ \ }K_{0}=3,K_{1}=1,K_{2}=3
\label{equati:pazertvbcunsmn}
\end{equation}%
respectively. Tribonacci concept was introduced by 14 year old student M.
Feinberg\ [\ref{bib:feinberg1963}] in 1963. Basic properties of it is given
in [\ref{bib:bruce1984}], [\ref{bib:scott1977}], [\ref{bib:shannon1977}], [%
\ref{bib:yalavigi1972}] and Binet formula for the $n$th number is given in [%
\ref{bib:spickerman1981}]. See also [\ref{catalani2002}],[\ref{bib:choi2013}%
],[\ref{bib:elia2001}], [\ref{bib:pethe1988}], [\ref{bib:lin1988}], [\ref%
{bib:spickerman1984}], [\ref{bib:yilmaz 2014}].

The sequences $\{T_{n}\}_{n\geq 0}$ and $\{K_{n}\}_{n\geq 0}$ can be
extended to negative subscripts by defining 
\begin{equation*}
T_{-n}=-T_{-(n-1)}-T_{-(n-2)}+T_{-(n-3)}
\end{equation*}%
and%
\begin{equation*}
K_{-n}=-K_{-(n-1)}-K_{-(n-2)}+K_{-(n-3)}
\end{equation*}%
for $n=1,2,3,...$ respectively. Therefore, recurrences (\ref%
{equati:fvcvxghsbnz}) and (\ref{equati:pazertvbcunsmn}) hold for all integer 
$n.$

It is well known that usual Tribonacci and Tribonacci-Lucas numbers can be
expressed using Binet's formulas%
\begin{equation}
T_{n}=\frac{\alpha ^{n+1}}{(\alpha -\beta )(\alpha -\gamma )}+\frac{\beta
^{n+1}}{(\beta -\alpha )(\beta -\gamma )}+\frac{\gamma ^{n+1}}{(\gamma
-\alpha )(\gamma -\beta )}  \label{equat:mnopcvbedcxzsa}
\end{equation}%
and%
\begin{equation*}
K_{n}=\alpha ^{n}+\beta ^{n}+\gamma ^{n}
\end{equation*}%
respectively, where $\alpha ,\beta $ and $\gamma $ are the roots of the
cubic equation $x^{3}-x^{2}-x-1=0.$ Moreover, 
\begin{eqnarray*}
\alpha &=&\frac{1+\sqrt[3]{19+3\sqrt{33}}+\sqrt[3]{19-3\sqrt{33}}}{3}, \\
\beta &=&\frac{1+\omega \sqrt[3]{19+3\sqrt{33}}+\omega ^{2}\sqrt[3]{19-3%
\sqrt{33}}}{3}, \\
\gamma &=&\frac{1+\omega ^{2}\sqrt[3]{19+3\sqrt{33}}+\omega \sqrt[3]{19-3%
\sqrt{33}}}{3}
\end{eqnarray*}%
where%
\begin{equation*}
\omega =\frac{-1+i\sqrt{3}}{2}=\exp (2\pi i/3),
\end{equation*}%
is a primitive cube root of unity. Note that we have the following identities%
\begin{eqnarray*}
\alpha +\beta +\gamma &=&1, \\
\alpha \beta +\alpha \gamma +\beta \gamma &=&-1, \\
\alpha \beta \gamma &=&1.
\end{eqnarray*}%
We can give Binet's formulas of the Tribonacci and Tribonacci-Lucas numbers
for the negative subscripts:\ for $n=1,2,3,...$ we have 
\begin{equation*}
T_{-n}=\frac{\alpha ^{-n+1}}{(\alpha -\beta )(\alpha -\gamma )}+\frac{\beta
^{-n+1}}{(\beta -\alpha )(\beta -\gamma )}+\frac{\gamma ^{-n+1}}{(\gamma
-\alpha )(\gamma -\beta )}
\end{equation*}%
and%
\begin{equation*}
K_{-n}=\alpha ^{-n}+\beta ^{-n}+\gamma ^{-n}.
\end{equation*}

The generating functions for the Tribonacci sequence $\{T_{n}\}_{n\geq 0}$
and Tribonacci-Lucas sequence $\{K_{n}\}_{n\geq 0}$ are 
\begin{equation*}
\sum_{n=0}^{\infty }T_{n}x^{n}=\frac{x}{1-x-x^{2}-x^{3}}\text{ \ and \ }%
\sum_{n=0}^{\infty }K_{n}x^{n}=\frac{3-2x-x^{2}}{1-x-x^{2}-x^{3}}.
\end{equation*}%
In this paper, we define Tribonacci and Tribonacci-Lucas sedenions in the
next section and give some properties of them. Before giving their
definition, we present some information on Cayley-Dickson algebras.

The algebras $\mathbb{C}$ (complex numbers), $\mathbb{H}$ (quaternions), and 
$\mathbb{O}$ (octonions) are real division algebras obtained from the real
numbers $\mathbb{R}$ by a doubling procedure called the Cayley-Dickson
Process (Construction). By doubling $\mathbb{R}$ (dim $2^{0}=1$), we obtain
the complex numbers $\mathbb{C}$ (dim $2^{1}=2$); then $\mathbb{C}$ yields
the quaternions $\mathbb{H}$ (dim $2^{2}=4$); and $\mathbb{H}$ produces
octonions $\mathbb{O}$ (dim $2^{3}=8$). The next doubling process applied to 
$\mathbb{O}$ then produces an algebra $\mathbb{S}$ (dim $2^{4}=16$) called
the sedenions. This doubling process can be extended beyond the sedenions to
form what are known as the $2^{n}$-ions (see\ for example [\ref%
{bib:imaeda2000}], [\ref{bib:moreno1998}], [\ref{bib:biss2008}]).

Next, we explain this doubling process.

The Cayley-Dickson algebras are a sequence $A_{0},A_{1},...$ of
non-associative $\mathbb{R}$-algebras with involution. The term
\textquotedblleft conjugation\textquotedblright\ can be used to refer to the
involution because it generalizes the usual conjugation on the complex
numbers. A full explanation of the basic properties of Cayley-Dickson
algebras, see [\ref{bib:biss2008}]. Cayley-Dickson algebras are defined
inductively. We begin by defining $A_{0}$ to be $\mathbb{R}.$ Given $%
A_{n-1}, $ the algebra $A_{n}$ is defined additively to be $A_{n-1}\times
A_{n-1}.$ Conjugation in $A_{n}$ is defined by 
\begin{equation*}
\overline{(a,b)}=(\overline{a},-b)
\end{equation*}%
and multiplication is defined by%
\begin{equation*}
(a,b)(c,d)=(ac-\overline{d}b,da+b\overline{c})
\end{equation*}%
and addition is defined by componentwise as%
\begin{equation*}
(a,b)+(c,d)=(a+c,b+d).
\end{equation*}%
Note that $A_{n}$ has dimension $2^{n}$ as an $\mathbb{R}-$vector space. If
we set, as usual, $\left\Vert x\right\Vert =\sqrt{\func{Re}(x\overline{x})}$
for $x\in A_{n}$ then $x\overline{x}=\overline{x}x=\left\Vert x\right\Vert
^{2}
.$

Now, suppose that $B_{16}=\{e_{i}\in \mathbb{S}:i=0,1,2,...,15\}$ is the
basis for $\mathbb{S}$, where $e_{0}$ is the identity (or unit) and $%
e_{1},e_{2},...,e_{15}$ are called imaginaries. Then a sedenion $S\in 
\mathbb{S}$ can be written as%
\begin{equation*}
S=\sum_{i=0}^{15}a_{i}e_{i}=a_{0}+\sum_{i=1}^{15}a_{i}e_{i}
\end{equation*}%
where $a_{0},a_{1},...,a_{15}$ are all real numbers. Here $a_{0}$ is called
the real part of $S$ and $\sum_{i=1}^{15}a_{i}e_{i}$ is called its imaginary
part.

Addition of sedenions is defined as componentwise and multiplication is
defined as follows: if $S_{1},S_{2}\in \mathbb{S}$ then we have%
\begin{equation}
S_{1}S_{2}=\left( \sum_{i=0}^{15}a_{i}e_{i}\right) \left(
\sum_{i=0}^{15}b_{i}e_{i}\right) =\sum_{i,j=0}^{15}a_{i}b_{j}(e_{i}e_{j}).
\label{equation:fvcxszeawqosbgv}
\end{equation}%
By setting $i\equiv e_{i}$ where $i=0,1,2,...,15,$ the multiplication rule
of the base elements $e_{i}\in B_{16}$ can be summarized as in the following
Table (see [\ref{bib:cawagas2004}] and [\ref{bib:bilgici2017}]). \FRAME{fhFU%
}{6.2699in}{2.3739in}{0pt}{\Qcb{Multiplication table for sedenions imaginary
units}}{}{tribonacci1.gif}{\special{language "Scientific Word";type
"GRAPHIC";maintain-aspect-ratio TRUE;display "USEDEF";valid_file "F";width
6.2699in;height 2.3739in;depth 0pt;original-width 6.2085in;original-height
2.3333in;cropleft "0";croptop "1";cropright "1";cropbottom "0";filename
'../tribonacci1.gif';file-properties "XNPEU";}}

From the above table, we can see that:

$e_{0}e_{i}=e_{i}e_{0}=e_{i};$ $e_{i}e_{i}=-e_{0}$ for $i\neq 0;$ $%
e_{i}e_{j}=-e_{j}e_{i}$ for $i\neq j$ and $i,j\neq 0.$

The operations requiring for the multiplication in (\ref%
{equation:fvcxszeawqosbgv}) are quite a lot. The computation of a sedenion
multiplication (product) using the naive method requires 256 multiplications
and 240 additions, while an algorithm which is given in [\ref{bib:cariow2013}%
] can compute the same result in only 122 multiplications (or multipliers --
in hardware implementation case) and 298 additions, for details see [\ref%
{bib:cariow2013}].

The problem with Cayley-Dickson Process is that each step of the doubling
process results in a progressive loss of structure. $\mathbb{R}$ is an
ordered field and it has all the nice properties we are so familiar with in
dealing with numbers like: the associative property, commutative property,
division property, self-conjugate property, etc. When we double $\mathbb{R}$
to have $\mathbb{C}$; $\mathbb{C}$ loses the self-conjugate property (and is
no longer an ordered field), next $\mathbb{H}$ loses the commutative
property, and $\mathbb{O}$ loses the associative property. When we double $%
\mathbb{O}$ to obtain $\mathbb{S}$; $\mathbb{S}$ loses the division
property. It means that $\mathbb{S}$ is non-commutative, non-associative,
and have a multiplicative identity element $e_{0}$ and multiplicative
inverses but it is not a division algebra because it has zero divisors; this
means that two non-zero sedenions can be multiplied to obtain zero: an
example is\ $(e_{3}+e_{10})(e_{6}-e_{15})=0$ and the other example is $%
(e_{2}-e_{14})(e_{3}+e_{15})=0,$ see [\ref{bib:cawagas2004}].

The algebras beyond the complex numbers go by the generic name hypercomplex
number. All hypercomplex number systems after sedenions that are based on
the Cayley--Dickson construction contain zero divisors.

Note that there is another type of sedenions which is called conic sedenions
or sedenions of Charles Muses, as they are also known, see [\ref%
{bib:koplinger2007a}], [\ref{bib:koplinger2007b}], [\ref{bib:muses1980}] for
more information. The term sedenion is also used for other 16-dimensional
algebraic structures, such as a tensor product of two copies of the
biquaternions, or the algebra of 4 by 4 matrices over the reals.

In the past, non-associative algebras and related structures with zero
divisors have not been given much attention because they did not appear to
have any useful applications in most mathematical subjects. Recently,
however, a lot of attention has been centred by theoretical physicists on
the Cayley-Dickson algebras $\mathbb{O}$ (octonions) and $\mathbb{S}$
(sedenions) because of their increasing usefulness in formulating many of
the new theories of elementary particles. In particular, the octonions $%
\mathbb{O}$ (which is the only non-associative normed division algebra over
the reals; see for example [\ref{bib:baez2002}] and [\ref{bib:okubo1995}])
has been found to be involved in so many unexpected areas (such as topology,
quantum theory, Clifford algebras, etc.) and sedenions appear in many areas
of science like linear gravity and electromagnetic theory.

Briefly $\mathbb{S}$, the algebra of sedenions, have the following
properties:

\begin{itemize}
\item $\mathbb{S}$ is a 16 dimensional non-associative and non-commutative
(Carley-Dickson) algebra over the reals,

\item $\mathbb{S}$ is not a composition algebra or division algebra because
of its zero divisors,

\item $\mathbb{S}$ is a non-alternative algebra, i.e., if $S_{1}$ and $S_{2}$
are sedenions the rules $S_{1}^{2}S_{2}=S_{1}(S_{1}S_{2})$ and $%
S_{1}S_{2}^{2}=(S_{1}S_{2})S_{2}$ do not always hold,

\item $\mathbb{S}$ is a power-associative algebra, i.e., if $S$ is an
sedenion then $S^{n}S^{m}=S^{n+m}.$
\end{itemize}

\section{The Tribonacci and Tribonacci-Lucas Sedenions and Their Generating
Functions and Binet Formulas}

In this section we define Tribonacci and Tribonacci-Lucas sedenions and give
generating functions and Binet formulas for them. First, we give some
information about quaternion sequences, octonion sequences and sedenion
sequences from literature.

Horadam [\ref{bib:horadam1963aa}] introduced $n$th Fibonacci and $n$th Lucas
quaternions as%
\begin{equation*}
Q_{n}=F_{n}+F_{n+1}e_{1}+F_{n+2}e_{2}+F_{n+3}e_{3}=\sum_{s=0}^{3}F_{n+s}e_{s}
\end{equation*}%
and%
\begin{equation*}
R_{n}=L_{n}+L_{n+1}e_{1}+L_{n+2}e_{2}+L_{n+3}e_{3}=\sum_{s=0}^{3}L_{n+s}e_{s}
\end{equation*}%
respectively, where $F_{n}$ and $L_{n}$ are the $n$th Fibonacci and Lucas
numbers respectively. He also defined generalized Fibonacci quaternion as%
\begin{equation*}
P_{n}=H_{n}+H_{n+1}e_{1}+H_{n+2}e_{2}+H_{n+3}e_{3}=\sum_{s=0}^{3}H_{n+s}e_{s}
\end{equation*}%
where $H_{n}$ is the $n$th generalized Fibonacci number (which is now called
Horadam number) by the recursive relation $H_{1}=p,$ $H_{2}=p+q,$ $%
H_{n}=H_{n-1}+H_{n-2}$ ($p$ and $q$ are arbitrary integers). Many other
generalization of Fibonacci quaternions has been given, see for example
Halici and Karata\c{s} [\ref{bib:halici2017}], and Polatl\i\ [\ref%
{bib:polatli2016}].

Cerda-Morales [\ref{bib:cerdamorale2017}] defined and studied the
generalized Tribonacci quaternion sequence that includes the previously
introduced Tribonacci, Padovan, Narayana and third order Jacobsthal
quaternion sequences. She defined generalized Tribonacci quaternion as%
\begin{equation*}
Q_{v,n}=V_{n}+V_{n+1}e_{1}+V_{n+2}e_{2}+V_{n+3}e_{3}=%
\sum_{s=0}^{3}V_{n+s}e_{s}
\end{equation*}%
where $V_{n}$ is the $n$th generalized Tribonacci number defined by the
third-order recurrance relations%
\begin{equation*}
V_{n}=rV_{n-1}+sV_{n-2}+tV_{n-3},\text{ \ \ \ \ }
\end{equation*}%
here $V_{0}=a,V_{1}=b,V_{2}=c$ are arbitrary integers and $r,s,t$ are real
numbers.

Various families of octonion number sequences (such as Fibonacci octonion,
Pell octonion, Jacobsthal octonion; and third order Jacobsthal octonion)
have been defined and studied by a number of authors in many different ways.
For example, Ke\c{c}ilioglu and Akku\c{s} [\ref{bib:kecilioglu}] introduced
the Fibonacci and Lucas octonions as%
\begin{equation*}
Q_{n}=\sum_{s=0}^{7}F_{n+s}e_{s}
\end{equation*}%
and%
\begin{equation*}
R_{n}=\sum_{s=0}^{7}L_{n+s}e_{s}
\end{equation*}%
respectively, where $F_{n}$ and $L_{n}$ are the $n$th Fibonacci and Lucas
numbers respectively. In [\ref{bib:cimen2017b}], \c{C}imen and \.{I}pek
introduced Jacobsthal octonions and Jacobsthal-Lucas octonions. In [\ref%
{bib:cerdamorale2017b}], Cerda-Morales introduced third order Jacobsthal
octonions and also in [\ref{bib:cerdamorale2018b}], she\ defined and studied
tribonacci-type octonions.

A number of authors have been defined and studied sedenion number sequences
(such as second order sedenions: Fibonacci sedenion, k-Pell and
k-Pell--Lucas sedenions, Jacobsthal and Jacobsthal-Lucas sedenions). For
example, Bilgici, Toke\c{s}er and \"{U}nal [\ref{bib:bilgici2017}]
introduced the Fibonacci and Lucas sedenions as%
\begin{equation*}
\widehat{F}_{n}=\sum_{s=0}^{15}F_{n+s}e_{s}
\end{equation*}%
and%
\begin{equation*}
\widehat{L}_{n}=\sum_{s=0}^{15}L_{n+s}e_{s}
\end{equation*}%
respectively, where $F_{n}$ and $L_{n}$ are the $n$th Fibonacci and Lucas
numbers respectively. In [\ref{bib:catarino2018}], Catarino introduced
k-Pell and k-Pell--Lucas sedenions. In [\ref{bib:cimen2017}], \c{C}imen and 
\.{I}pek introduced Jacobsthal and Jacobsthal-Lucas sedenions.

G\"{u}l [\ref{bib:gul2018}] introduced the k-Fibonacci and k-Lucas
trigintaduonions as%
\begin{equation*}
TF_{k,n}=\sum_{s=0}^{31}F_{k,n+s}e_{s}
\end{equation*}%
and%
\begin{equation*}
TL_{k,n}=\sum_{s=0}^{31}L_{k,n+s}e_{s}
\end{equation*}%
respectively, where $F_{k,n}$ and $L_{k,n}$ are the $n$th k-Fibonacci and
k-Lucas numbers respectively.

We now define Tribonacci and Tribonacci-Lucas sedenions over the sedenion
algebra $\mathbb{S}$. The $n$th Tribonacci sedenion is%
\begin{equation}
\widehat{T}_{n}=\sum_{s=0}^{15}T_{n+s}e_{s}=T_{n}+\sum_{s=1}^{15}T_{n+s}e_{s}
\label{equation:bnuyvbcxdszaerm}
\end{equation}%
and the $n$th Tribonacci-Lucas sedenion is%
\begin{equation*}
\widehat{K}_{n}=\sum_{s=0}^{15}K_{n+s}e_{s}=K_{n}+%
\sum_{s=1}^{15}K_{n+s}e_{s}.
\end{equation*}%
For all non-negative integer $n,$ it can be easily shown that 
\begin{equation}
\widehat{T}_{n}=\widehat{T}_{n-1}+\widehat{T}_{n-2}+\widehat{T}_{n-3}
\label{equation:abgsvbnuytsaerds}
\end{equation}%
and%
\begin{equation}
\widehat{K}_{n}=\widehat{K}_{n-1}+\widehat{K}_{n-2}+\widehat{K}_{n-3}.
\label{equation:cvnbhustarfcxdsa}
\end{equation}
The sequences $\{\widehat{T}_{n}\}_{n\geq 0}$ and $\{\widehat{K}%
_{n}\}_{n\geq 0}$ can be defined for negative values of $n$ by using the
recurrences (\ref{equation:abgsvbnuytsaerds}) and (\ref%
{equation:cvnbhustarfcxdsa}) to extend the sequence backwards, or
equivalently, by using the recurrences 
\begin{equation*}
\widehat{T}_{-n}=-\widehat{T}_{-(n-1)}-\widehat{T}_{-(n-2)}+\widehat{T}%
_{-(n-3)}
\end{equation*}%
and%
\begin{equation*}
\widehat{K}_{-n}=-\widehat{K}_{-(n-1)}-\widehat{K}_{-(n-2)}+\widehat{K}%
_{-(n-3)},
\end{equation*}%
respectively. Thus, the recurrences (\ref{equation:abgsvbnuytsaerds}) and (%
\ref{equation:cvnbhustarfcxdsa}) holds for all integer $n.$

The conjugate of $\widehat{T}_{n}$ and $\widehat{K}_{n}$ are defined by%
\begin{equation*}
\overline{\widehat{T}_{n}}=T_{n}-%
\sum_{s=1}^{15}T_{n+s}e_{s}=T_{n}-T_{n+1}e_{1}-T_{n+2}e_{2}-...-T_{n+15}e_{15}
\end{equation*}%
and%
\begin{equation*}
\overline{\widehat{K}_{n}}=K_{n}-%
\sum_{s=1}^{15}K_{n+s}e_{s}=K_{n}-K_{n+1}e_{1}-K_{n+2}e_{2}-...-K_{n+15}e_{15}
\end{equation*}%
respectively. The norms of $n$th Tribonacci and Tribonacci-Lucas sedenions
are 
\begin{equation*}
\left\Vert \widehat{T}_{n}\right\Vert ^{2}=N^{2}(\widehat{T}_{n})=\widehat{T}%
_{n}\overline{\widehat{T}_{n}}=\overline{\widehat{T}_{n}}\widehat{T}%
_{n}=T_{n}^{2}+T_{n+1}^{2}+...+T_{n+15}^{2}
\end{equation*}%
and%
\begin{equation*}
\left\Vert \widehat{K}_{n}\right\Vert ^{2}=N^{2}(\widehat{K}_{n})=\widehat{K}%
_{n}\overline{\widehat{K}_{n}}=\overline{\widehat{K}_{n}}\widehat{K}%
_{n}=K_{n}^{2}+K_{n+1}^{2}+...+K_{n+15}^{2}
\end{equation*}%
respectively.

Now, we will state Binet's formula for the Tribonacci and Tribonacci-Lucas
sedenions and in the rest of the paper we fixed the following notations.%
\begin{eqnarray*}
\widehat{\alpha } &=&\sum_{s=0}^{15}\alpha ^{s}e_{s}, \\
\widehat{\beta } &=&\sum_{s=0}^{15}\beta ^{s}e_{s}, \\
\widehat{\gamma } &=&\sum_{s=0}^{15}\gamma ^{s}e_{s}.
\end{eqnarray*}

\begin{theorem}
\label{theorem:yunhbvgcfosea}For any integer $n,$ the $n$th Tribonacci
sedenion is%
\begin{equation}
\widehat{T}_{n}=\frac{\widehat{\alpha }\alpha ^{n+1}}{(\alpha -\beta
)(\alpha -\gamma )}+\frac{\widehat{\beta }\beta ^{n+1}}{(\beta -\alpha
)(\beta -\gamma )}+\frac{\widehat{\gamma }\gamma ^{n+1}}{(\gamma -\alpha
)(\gamma -\beta )}  \label{equation:fgvbvcxsdaewqzxd}
\end{equation}%
and the $n$th Tribonacci-Lucas sedenion is 
\begin{equation}
\widehat{K}_{n}=\widehat{\alpha }\alpha ^{n}+\widehat{\beta }\beta ^{n}+%
\widehat{\gamma }\gamma ^{n}.  \label{equation:xcdfvesaztfvop}
\end{equation}
\end{theorem}

\textit{Proof.} Repeated use of (\ref{equat:mnopcvbedcxzsa}) in (\ref%
{equation:bnuyvbcxdszaerm}) enable us to write for $\widehat{\alpha }%
=\sum_{s=0}^{15}\alpha ^{s}e_{s},$ $\widehat{\beta }=\sum_{s=0}^{15}\beta
^{s}e_{s}$ and $\widehat{\gamma }=\sum_{s=0}^{15}\gamma ^{s}e_{s}:$%
\begin{eqnarray*}
\widehat{T}_{n} &=&\sum_{s=0}^{15}T_{n+s}e_{s}=\sum_{s=0}^{15}\left( \frac{%
\alpha ^{n+1+s}e_{s}}{(\alpha -\beta )(\alpha -\gamma )}+\frac{\beta
^{n+1+s}e_{s}}{(\beta -\alpha )(\beta -\gamma )}+\frac{\gamma ^{n+1+s}e_{s}}{%
(\gamma -\alpha )(\gamma -\beta )}\right) \\
&=&\frac{\widehat{\alpha }\alpha ^{n+1}}{(\alpha -\beta )(\alpha -\gamma )}+%
\frac{\widehat{\beta }\beta ^{n+1}}{(\beta -\alpha )(\beta -\gamma )}+\frac{%
\widehat{\gamma }\gamma ^{n+1}}{(\gamma -\alpha )(\gamma -\beta )}.
\end{eqnarray*}%
Similarly, we can obtain (\ref{equation:xcdfvesaztfvop}). 
\endproof%

We can give Binet's formula of the Tribonacci and Tribonacci-Lucas sedenions
for the negative subscripts:\ for $n=1,2,3,...$ we have%
\begin{equation*}
\widehat{T}_{-n}=\frac{\widehat{\alpha }\alpha ^{-n+1}}{(\alpha -\beta
)(\alpha -\gamma )}+\frac{\widehat{\beta }\beta ^{-n+1}}{(\beta -\alpha
)(\beta -\gamma )}+\frac{\widehat{\gamma }\gamma ^{-n+1}}{(\gamma -\alpha
)(\gamma -\beta )}
\end{equation*}%
and%
\begin{equation*}
\widehat{K}_{-n}=\widehat{\alpha }\alpha ^{-n}+\widehat{\beta }\beta ^{-n}+%
\widehat{\gamma }\gamma ^{-n},
\end{equation*}%
respectively.

The next theorem gives us an alternatif proof of the Binet's formula for the
Tribonacci and Tribonacci-Lucas sedenions. For this, we need the quadratik
approximation of $\{T_{n}\}_{n\geq 0}$ and $\{K_{n}\}_{n\geq 0}$:

\begin{lemma}
For all integer $n,$ we have

\begin{description}
\item[(a)] 
\begin{eqnarray*}
\alpha \alpha ^{n+2} &=&T_{n+2}\alpha ^{2}+(T_{n+1}+T_{n})\alpha +T_{n+1}, \\
\beta \beta ^{n+2} &=&T_{n+2}\beta ^{2}+(T_{n+1}+T_{n})\beta +T_{n+1}, \\
\gamma \gamma ^{n+2} &=&T_{n+2}\gamma ^{2}+(T_{n+1}+T_{n})\gamma +T_{n+1}.
\end{eqnarray*}

\item[(b)] 
\begin{eqnarray*}
P\alpha ^{n+2} &=&K_{n+2}\alpha ^{2}+(K_{n+1}+K_{n})\alpha +K_{n+1}, \\
Q\beta ^{n+2} &=&K_{n+2}\beta ^{2}+(K_{n+1}+K_{n})\beta +K_{n+1}, \\
R\gamma ^{n+2} &=&K_{n+2}\gamma ^{2}+(K_{n+1}+K_{n})\gamma +K_{n+1},
\end{eqnarray*}%
where%
\begin{eqnarray*}
P &=&3-(\beta +\gamma )+3\beta \gamma , \\
Q &=&3-(\alpha +\gamma )+3\alpha \gamma , \\
R &=&3-(\alpha +\beta \gamma )+3\alpha \beta .
\end{eqnarray*}
\end{description}
\end{lemma}

\textit{Proof.} See [\ref{bib:cerdamorale2018a}] or [\ref%
{bib:cerdamorale2018b}]. 
\endproof%

\underline{Alternatif Proof of Theorem \ref{theorem:yunhbvgcfosea}:}

Note that%
\begin{eqnarray*}
&&\alpha ^{2}\widehat{T}_{n+2}+\alpha (\widehat{T}_{n+1}+\widehat{T}_{n})+%
\widehat{T}_{n+1} \\
&=&\alpha ^{2}(T_{n+2}+T_{n+3}e_{1}+...+T_{n+17}e_{15}) \\
&&+\alpha
((T_{n+1}+T_{n})+(T_{n+2}+T_{n+1})e_{1}+...+(T_{n+16}+T_{n+15})e_{15}) \\
&&+(T_{n+1}+T_{n+2}e_{1}+...+T_{n+16}e_{15}) \\
&=&\alpha ^{2}T_{n+2}+\alpha (T_{n+1}+T_{n})+T_{n+1}+(\alpha
^{2}T_{n+3}+(T_{n+2}+T_{n+1})+T_{n+2})e_{1} \\
&&+(\alpha ^{2}T_{n+4}+(T_{n+3}+T_{n+2})+T_{n+3})e_{2} \\
&&\vdots \\
&&+(\alpha ^{2}T_{n+17}+(T_{n+16}+T_{n+15})+T_{n+16})e_{15}.
\end{eqnarray*}%
From the identity $\alpha ^{n+3}=T_{n+2}\alpha ^{2}+(T_{n+1}+T_{n})\alpha
+T_{n+1}$ for $n$-th Tribonacci number $T_{n},$ we have%
\begin{equation}
\alpha ^{2}\widehat{T}_{n+2}+\alpha (\widehat{T}_{n+1}+\widehat{T}_{n})+%
\widehat{T}_{n+1}=\widehat{\alpha }\alpha ^{n+3}.
\label{equation:abctyunmbvcfoes}
\end{equation}%
Similarly, we obtain%
\begin{equation}
\beta ^{2}\widehat{T}_{n+2}+\beta (\widehat{T}_{n+1}+\widehat{T}_{n})+%
\widehat{T}_{n+1}=\widehat{\beta }\beta ^{n+3}
\label{equation:byummnbvoepsdxcz}
\end{equation}%
and%
\begin{equation}
\gamma ^{2}\widehat{T}_{n+2}+\gamma (\widehat{T}_{n+1}+\widehat{T}_{n})+%
\widehat{T}_{n+1}=\widehat{\gamma }\gamma ^{n+3}.
\label{equation:cuytnbhvgfopsexza}
\end{equation}%
Substracting (\ref{equation:byummnbvoepsdxcz}) from (\ref%
{equation:abctyunmbvcfoes}), we have%
\begin{equation}
(\alpha +\beta )\widehat{T}_{n+2}+(\widehat{T}_{n+1}+\widehat{T}_{n})=\frac{%
\widehat{\alpha }\alpha ^{n+3}-\widehat{\beta }\beta ^{n+3}}{\alpha -\beta }.
\label{equation:somnbyusxczaeo}
\end{equation}%
Similarly, substracting (\ref{equation:cuytnbhvgfopsexza}) from (\ref%
{equation:abctyunmbvcfoes}), we obtain%
\begin{equation}
(\alpha +\gamma )\widehat{T}_{n+2}+(\widehat{T}_{n+1}+\widehat{T}_{n})=\frac{%
\widehat{\alpha }\alpha ^{n+3}-\widehat{\gamma }\gamma ^{n+3}}{\alpha
-\gamma }.  \label{equation:gtybnvcxdsxzc}
\end{equation}%
Finally, substracting (\ref{equation:gtybnvcxdsxzc}) from (\ref%
{equation:somnbyusxczaeo}), we get%
\begin{eqnarray*}
\widehat{T}_{n+2} &=&\frac{1}{\alpha -\beta }\left( \frac{\widehat{\alpha }%
\alpha ^{n+3}-\widehat{\beta }\beta ^{n+3}}{\alpha -\beta }-\frac{\widehat{%
\alpha }\alpha ^{n+3}-\widehat{\gamma }\gamma ^{n+3}}{\alpha -\gamma }\right)
\\
&=&\frac{\widehat{\alpha }\alpha ^{n+3}}{(\alpha -\beta )(\alpha -\gamma )}-%
\frac{\widehat{\beta }\beta ^{n+3}}{(\alpha -\beta )(\beta -\gamma )}+\frac{%
\widehat{\gamma }\gamma ^{n+3}}{(\gamma -\alpha )(\gamma -\beta )} \\
&=&\frac{\widehat{\alpha }\alpha ^{n+3}}{(\alpha -\beta )(\alpha -\gamma )}+%
\frac{\widehat{\beta }\beta ^{n+3}}{(\beta -\alpha )(\beta -\gamma )}+\frac{%
\widehat{\gamma }\gamma ^{n+3}}{(\gamma -\alpha )(\gamma -\beta )}.
\end{eqnarray*}%
So this proves (\ref{equation:fgvbvcxsdaewqzxd}). Similarly we obtain (\ref%
{equation:xcdfvesaztfvop}). 
\endproof%

Next, we present generating functions.

\begin{theorem}
The generating functions for the Tribonacci and Tribonacci-Lucas sedenions
are%
\begin{equation}
g(x)=\sum_{n=0}^{\infty }\widehat{T}_{n}x^{n}=\frac{\widehat{T}_{0}+(%
\widehat{T}_{1}-\widehat{T}_{0})x+\widehat{T}_{-1}x^{2}}{1-x-x^{2}-x^{3}}
\label{equation:hnbvyuedsxzapo}
\end{equation}%
and%
\begin{equation}
g(x)=\sum_{n=0}^{\infty }\widehat{K}_{n}x^{n}=\frac{\widehat{K}_{0}+(%
\widehat{K}_{1}-\widehat{K}_{0})x+\widehat{K}_{-1}x^{2}}{1-x-x^{2}-x^{3}}
\label{equat:bnopscxzartfnm}
\end{equation}%
respectively.
\end{theorem}

\textit{Proof.} Define $g(x)=\sum_{n=0}^{\infty }\widehat{T}_{n}x^{n}.$ Note
that%
\begin{equation*}
\begin{array}{ccccccccc}
g(x) & = & \widehat{T}_{0}+ & \widehat{T}_{1}x+ & \widehat{T}_{2}x^{2}+ & 
\widehat{T}_{3}x^{3}+ & \widehat{T}_{4}x^{4}+ & \widehat{T}_{5}x^{5}+...+ & 
\widehat{T}_{n}x^{n}\text{ }+... \\ 
xg(x) & = &  & \widehat{T}_{0}x+ & \widehat{T}_{1}x^{2}+ & \widehat{T}%
_{2}x^{3}+ & \widehat{T}_{3}x^{4}+ & \widehat{T}_{4}x^{5}+...+ & \widehat{T}%
_{n-1}x^{n}+... \\ 
x^{2}g(x) & = &  &  & \widehat{T}_{0}x^{2}+ & \widehat{T}_{1}x^{3}+ & 
\widehat{T}_{2}x^{4}+ & \widehat{T}_{3}x^{5}+...+ & \widehat{T}%
_{n-2}x^{n}+... \\ 
x^{3}g(x) & = &  &  &  & \widehat{T}_{0}x^{3}+ & \widehat{T}_{1}x^{4}+ & 
\widehat{T}_{2}x^{5}+...+ & \widehat{T}_{n-3}x^{n}+...%
\end{array}%
\end{equation*}%
Using above table and the recurans $\widehat{T}_{n}=\widehat{T}_{n-1}+%
\widehat{T}_{n-2}+\widehat{T}_{n-3}$ we have 
\begin{eqnarray*}
&&g(x)-xg(x)-x^{2}g(x)-x^{3}g(x) \\
&=&\widehat{T}_{0}+(\widehat{T}_{1}-\widehat{T}_{0})x+(\widehat{T}_{2}-%
\widehat{T}_{1}-\widehat{T}_{0})x^{2}+(\widehat{T}_{3}-\widehat{T}_{2}-%
\widehat{T}_{1}-\widehat{T}_{0})x^{3}+ \\
&&(\widehat{T}_{4}-\widehat{T}_{3}-\widehat{T}_{2}-\widehat{T}%
_{1})x^{4}+...+(\widehat{T}_{n}-\widehat{T}_{n-1}-\widehat{T}_{n-2}-\widehat{%
T}_{n-3}+)x^{n}+... \\
&=&\widehat{T}_{0}+(\widehat{T}_{1}-\widehat{T}_{0})x+(\widehat{T}_{2}-%
\widehat{T}_{1}-\widehat{T}_{0})x^{2}.
\end{eqnarray*}%
It follows that%
\begin{equation*}
g(x)=\frac{\widehat{T}_{0}+(\widehat{T}_{1}-\widehat{T}_{0})x+(\widehat{T}%
_{2}-\widehat{T}_{1}-\widehat{T}_{0})x^{2}}{1-x-x^{2}-x^{3}}.
\end{equation*}%
Since $\widehat{T}_{2}-\widehat{T}_{1}-\widehat{T}_{0}=\widehat{T}_{-1},$
the generating functions for the Tribonacci sedenion is%
\begin{equation*}
g(x)=\frac{\widehat{T}_{0}+(\widehat{T}_{1}-\widehat{T}_{0})x+\widehat{T}%
_{-1}x^{2}}{1-x-x^{2}-x^{3}}.
\end{equation*}%
Similarly, we can obtain (\ref{equation:hnbvyuedsxzapo}). 
\endproof%

In the following theorem we present another forms of Binet formulas for the
Tribonacci and Tribonacci-Lucas sedenions using generating functions.

\begin{theorem}
\label{theorem:vbghcxdszrtupomn}For any integer $n,$ the $n$th Tribonacci
sedenion is%
\begin{eqnarray*}
\widehat{T}_{n} &=&\frac{((\alpha ^{2}-\alpha )\widehat{T}_{0}+\alpha 
\widehat{T}_{1}+\widehat{T}_{-1})\alpha ^{n}}{\left( \alpha -\gamma \right)
\left( \alpha -\beta \right) }+\frac{((\beta ^{2}-\beta )\widehat{T}%
_{0}+\beta \widehat{T}_{1}+\widehat{T}_{-1})\beta ^{n}}{\left( \beta -\gamma
\right) \left( \beta -\alpha \right) } \\
&&+\frac{((\gamma ^{2}-\gamma )\widehat{T}_{0}+\gamma \widehat{T}_{1}+%
\widehat{T}_{-1})\gamma ^{n}}{\left( \gamma -\beta \right) \left( \gamma
-\alpha \right) }
\end{eqnarray*}%
and the $n$th Tribonacci-Lucas sedenion is%
\begin{eqnarray*}
\widehat{K}_{n} &=&\frac{((\alpha ^{2}-\alpha )\widehat{K}_{0}+\alpha 
\widehat{K}_{1}+\widehat{K}_{-1})\alpha ^{n}}{\left( \alpha -\gamma \right)
\left( \alpha -\beta \right) }+\frac{((\beta ^{2}-\beta )\widehat{K}%
_{0}+\beta \widehat{K}_{1}+\widehat{K}_{-1})\beta ^{n}}{\left( \beta -\gamma
\right) \left( \beta -\alpha \right) } \\
&&+\frac{((\gamma ^{2}-\gamma )\widehat{K}_{0}+\gamma \widehat{K}_{1}+%
\widehat{K}_{-1})\gamma ^{n}}{\left( \gamma -\beta \right) \left( \gamma
-\alpha \right) }.
\end{eqnarray*}
\end{theorem}

\textit{Proof.} We can use generating functions. Since the roots of the
equation $1-x-x^{2}-x^{3}=0$ are $\alpha \beta ,\beta \gamma ,\alpha \gamma $
and 
\begin{equation*}
1-x-x^{2}-x^{3}=(1-\alpha x)(1-\beta x)(1-\gamma x)
\end{equation*}%
we can write the generating function of $\widehat{T}_{n}$ as \ 
\begin{eqnarray*}
&&g(x) \\
&=&\frac{\widehat{T}_{0}+(\widehat{T}_{1}-\widehat{T}_{0})x+\widehat{T}%
_{-1}x^{2}}{1-x-x^{2}-x^{3}}=\frac{\widehat{T}_{0}+(\widehat{T}_{1}-\widehat{%
T}_{0})x+\widehat{T}_{-1}x^{2}}{(1-\alpha x)(1-\beta x)(1-\gamma x)} \\
&=&\frac{A}{(1-\alpha x)}+\frac{B}{(1-\beta x)}+\frac{C}{(1-\gamma x)} \\
&=&\frac{A(1-\beta x)(1-\gamma x)+B(1-\alpha x)(1-\gamma x)+C(1-\alpha
x)(1-\beta x)}{(1-\alpha x)(1-\beta x)(1-\gamma x)} \\
&=&\frac{(A+B+C)+(-A\beta -A\gamma -B\alpha -B\gamma -C\alpha -C\beta
)x+(A\beta \gamma +B\alpha \gamma +C\alpha \beta )x^{2}}{(1-\alpha
x)(1-\beta x)(1-\gamma x)}.
\end{eqnarray*}%
We need to find $A,B$ and $C$, so the following system of equations should
be solved:%
\begin{eqnarray*}
A+B+C &=&\widehat{T}_{0} \\
-A\beta -A\gamma -B\alpha -B\gamma -C\alpha -C\beta &=&\widehat{T}_{1}-%
\widehat{T}_{0} \\
A\beta \gamma +B\alpha \gamma +C\alpha \beta &=&\widehat{T}_{-1}
\end{eqnarray*}%
We find that%
\begin{eqnarray*}
A &=&\frac{-\alpha \widehat{T}_{0}+\alpha \widehat{T}_{1}+\widehat{T}%
_{-1}+\alpha ^{2}\widehat{T}_{0}}{\alpha ^{2}-\alpha \beta -\alpha \gamma
+\beta \gamma }=\frac{((\alpha ^{2}-\alpha )\widehat{T}_{0}+\alpha \widehat{T%
}_{1}+\widehat{T}_{-1})}{\left( \alpha -\gamma \right) \left( \alpha -\beta
\right) }, \\
B &=&\frac{-\beta \widehat{T}_{0}+\beta \widehat{T}_{1}+\widehat{T}%
_{-1}+\beta ^{2}\widehat{T}_{0}}{\beta ^{2}-\alpha \beta +\alpha \gamma
-\beta \gamma }=\frac{((\beta ^{2}-\beta )\widehat{T}_{0}+\beta \widehat{T}%
_{1}+\widehat{T}_{-1})}{\left( \beta -\gamma \right) \left( \beta -\alpha
\right) }, \\
C &=&\frac{-\gamma \widehat{T}_{0}+\gamma \widehat{T}_{1}+\widehat{T}%
_{-1}+\gamma ^{2}\widehat{T}_{0}}{\gamma ^{2}+\alpha \beta -\alpha \gamma
-\beta \gamma }=\frac{((\gamma ^{2}-\gamma )\widehat{T}_{0}+\gamma \widehat{T%
}_{1}+\widehat{T}_{-1})}{\left( \gamma -\beta \right) \left( \gamma -\alpha
\right) }.
\end{eqnarray*}%
and 
\begin{eqnarray*}
g(x) &=&\frac{((\alpha ^{2}-\alpha )\widehat{T}_{0}+\alpha \widehat{T}_{1}+%
\widehat{T}_{-1})}{\left( \alpha -\gamma \right) \left( \alpha -\beta
\right) }\sum_{n=0}^{\infty }\alpha ^{n}x^{n}+\frac{((\beta ^{2}-\beta )%
\widehat{T}_{0}+\beta \widehat{T}_{1}+\widehat{T}_{-1})}{\left( \beta
-\gamma \right) \left( \beta -\alpha \right) }\sum_{n=0}^{\infty }\beta
^{n}x^{n} \\
&&+\frac{(-\gamma \widehat{T}_{0}+\gamma \widehat{T}_{1}+\widehat{T}%
_{-1}+\gamma ^{2}\widehat{T}_{0})}{\left( \gamma -\beta \right) \left(
\gamma -\alpha \right) }\sum_{n=0}^{\infty }\gamma ^{n}x^{n} \\
&=&\sum_{n=0}^{\infty }\left( 
\begin{array}{c}
\frac{((\alpha ^{2}-\alpha )\widehat{T}_{0}+\alpha \widehat{T}_{1}+\widehat{T%
}_{-1})\alpha ^{n}}{\left( \alpha -\gamma \right) \left( \alpha -\beta
\right) }+\frac{((\beta ^{2}-\beta )\widehat{T}_{0}+\beta \widehat{T}_{1}+%
\widehat{T}_{-1})\beta ^{n}}{\left( \beta -\gamma \right) \left( \beta
-\alpha \right) } \\ 
+\frac{((\gamma ^{2}-\gamma )\widehat{T}_{0}+\gamma \widehat{T}_{1}+\widehat{%
T}_{-1})\gamma ^{n}}{\left( \gamma -\beta \right) \left( \gamma -\alpha
\right) }%
\end{array}%
\right) x^{n}.
\end{eqnarray*}%
Thus Binet formula of Tribonacci sedenion is 
\begin{eqnarray*}
\widehat{T}_{n} &=&\frac{((\alpha ^{2}-\alpha )\widehat{T}_{0}+\alpha 
\widehat{T}_{1}+\widehat{T}_{-1})\alpha ^{n}}{\left( \alpha -\gamma \right)
\left( \alpha -\beta \right) }+\frac{((\beta ^{2}-\beta )\widehat{T}%
_{0}+\beta \widehat{T}_{1}+\widehat{T}_{-1})\beta ^{n}}{\left( \beta -\gamma
\right) \left( \beta -\alpha \right) } \\
&&+\frac{((\gamma ^{2}-\gamma )\widehat{T}_{0}+\gamma \widehat{T}_{1}+%
\widehat{T}_{-1})\gamma ^{n}}{\left( \gamma -\beta \right) \left( \gamma
-\alpha \right) }.
\end{eqnarray*}%
Similarly, we can obtain Binet formula of the Tribonacci-Lucas sedenion. 
\endproof%

If we compare Theorem \ref{theorem:yunhbvgcfosea} and Theorem \ref%
{theorem:vbghcxdszrtupomn} and use the definition of $\widehat{T}_{n},$ $%
\widehat{K}_{n},$ we have the following Remark showing relations between $%
\widehat{T}_{-1},\widehat{T}_{0},\widehat{T}_{1};\widehat{K}_{-1},\widehat{K}%
_{0},\widehat{K}_{1}$ and $\widehat{\alpha },\widehat{\beta },\widehat{%
\gamma }.$ We obtain (b) and (d) after solving the system of the equations
in (a) and (b) respectively.

\begin{remark}
We have the following identities:

\begin{description}
\item[(a)] 
\begin{eqnarray*}
\frac{(\alpha ^{2}-\alpha )\widehat{T}_{0}+\alpha \widehat{T}_{1}+\widehat{T}%
_{-1}}{\alpha } &=&\widehat{\alpha } \\
\frac{(\beta ^{2}-\beta )\widehat{T}_{0}+\beta \widehat{T}_{1}+\widehat{T}%
_{-1}}{\beta } &=&\widehat{\beta } \\
\frac{(\gamma ^{2}-\gamma )\widehat{T}_{0}+\gamma \widehat{T}_{1}+\widehat{T}%
_{-1}}{\gamma } &=&\widehat{\gamma }
\end{eqnarray*}

\item[(b)] 
\begin{eqnarray*}
\sum_{s=0}^{15}T_{-1+s}e_{s} &=&\widehat{T}_{-1}=\frac{\widehat{\alpha }}{%
(\alpha -\beta )(\alpha -\gamma )}+\frac{\widehat{\beta }}{(\beta -\alpha
)(\beta -\gamma )}+\frac{\widehat{\gamma }}{(\gamma -\alpha )(\gamma -\beta )%
} \\
\sum_{s=0}^{15}T_{s}e_{s} &=&\widehat{T}_{0}=\dfrac{\widehat{\alpha }\alpha 
}{(\alpha -\beta )(\alpha -\gamma )}+\dfrac{\widehat{\beta }\beta }{(\beta
-\alpha )(\beta -\gamma )}+\dfrac{\widehat{\gamma }\gamma }{(\gamma -\alpha
)(\gamma -\beta )} \\
\sum_{s=0}^{15}T_{1+s}e_{s} &=&\widehat{T}_{1}=\frac{\widehat{\alpha }\alpha
^{2}}{(\alpha -\beta )(\alpha -\gamma )}+\frac{\widehat{\beta }\beta ^{2}}{%
(\beta -\alpha )(\beta -\gamma )}+\frac{\widehat{\gamma }\gamma ^{2}}{%
(\gamma -\alpha )(\gamma -\beta )}
\end{eqnarray*}

\item[(c)] 
\begin{eqnarray*}
\dfrac{((\alpha ^{2}-\alpha )\widehat{K}_{0}+\alpha \widehat{K}_{1}+\widehat{%
K}_{-1})}{\left( \alpha -\gamma \right) \left( \alpha -\beta \right) } &=&%
\widehat{\alpha } \\
\dfrac{((\beta ^{2}-\beta )\widehat{K}_{0}+\beta \widehat{K}_{1}+\widehat{K}%
_{-1})}{\left( \beta -\gamma \right) \left( \beta -\alpha \right) } &=&%
\widehat{\beta } \\
\dfrac{((\gamma ^{2}-\gamma )\widehat{K}_{0}+\gamma \widehat{K}_{1}+\widehat{%
K}_{-1})}{\left( \gamma -\beta \right) \left( \gamma -\alpha \right) } &=&%
\widehat{\gamma }
\end{eqnarray*}

\item[(d)] 
\begin{eqnarray*}
\sum_{s=0}^{15}K_{-1+s}e_{s} &=&\widehat{K}_{-1}=\widehat{\alpha }\alpha
^{-1}+\widehat{\beta }\beta ^{-1}+\widehat{\gamma }\gamma ^{-1} \\
\sum_{s=0}^{15}K_{s}e_{s} &=&\widehat{K}_{0}=\widehat{\alpha }+\widehat{%
\beta }+\widehat{\gamma } \\
\sum_{s=0}^{15}K_{1+s}e_{s} &=&\widehat{K}_{1}=\widehat{\alpha }\alpha +%
\widehat{\beta }\beta +\widehat{\gamma }\gamma .
\end{eqnarray*}
\end{description}
\end{remark}

Using above Remark we can find $\widehat{T}_{2},$ $\widehat{K}_{2}$ as
follows:%
\begin{equation}
\sum_{s=0}^{15}T_{2+s}e_{s}=\widehat{T}_{2}=\widehat{T}_{1}+\widehat{T}_{0}+%
\widehat{T}_{-1}=\frac{\widehat{\alpha }\alpha ^{3}}{(\alpha -\beta )(\alpha
-\gamma )}+\frac{\widehat{\beta }\beta ^{3}}{(\beta -\alpha )(\beta -\gamma )%
}+\frac{\widehat{\gamma }\gamma ^{3}}{(\gamma -\alpha )(\gamma -\beta )}
\label{equati:mnosratvbxczgus}
\end{equation}

and%
\begin{equation}
\sum_{s=0}^{15}K_{2+s}e_{s}=\widehat{K}_{2}=\widehat{K}_{1}+\widehat{K}_{0}+%
\widehat{K}_{-1}=\widehat{\alpha }\alpha ^{2}+\widehat{\beta }\beta ^{2}+%
\widehat{\gamma }\gamma ^{2}.  \label{equati:mnbvyughbqazse}
\end{equation}

Of course, (\ref{equati:mnosratvbxczgus}) and (\ref{equati:mnbvyughbqazse})
can be found directly from (\ref{equation:fgvbvcxsdaewqzxd}) and (\ref%
{equation:xcdfvesaztfvop}).

Now, we present the formulas which give the summation of the first $n$
Tribonacci and Tribonacci-Lucas numbers.

\begin{lemma}
For every integer $n\geq 0,$ we have 
\begin{equation}
\sum\limits_{i=0}^{n}T_{i}=T_{0}+\dfrac{1}{2}(T_{n+2}+T_{n}-1)=\dfrac{1}{2}%
(T_{n+2}+T_{n}-1)  \label{equation:easzpuyhbnvcdf}
\end{equation}%
and%
\begin{equation}
\sum\limits_{i=0}^{n}K_{i}=\dfrac{K_{n+2}+K_{n}}{2}.
\label{equation:prtsaeoaxcdf}
\end{equation}
\end{lemma}

\textit{Proof.} (\ref{equation:easzpuyhbnvcdf}) and (\ref%
{equation:prtsaeoaxcdf}) can be proved by mathematical induction easily. For
a proof of (\ref{equation:easzpuyhbnvcdf}) with a telescopik sum method see [%
\ref{bib:devbhadra2011}] or with a matrix diagonalisation proof, see [\ref%
{bib:kilic2008}] or see also [\ref{bib:cerdamorale2017}].

For a proof of (\ref{equation:prtsaeoaxcdf}), see [\ref{bib:frontczak2018}].
Since $K_{0}=3$ and $\sum\limits_{i=1}^{n}K_{i}=\dfrac{K_{n+2}+K_{n}-6}{2},$
it follows that $\sum\limits_{i=0}^{n}K_{i}=\dfrac{K_{n+2}+K_{n}}{2}.$ 
\endproof%

There is also a formula of the summation of the first $n$ negative
Tribonacci numbers:%
\begin{equation*}
\sum\limits_{i=1}^{n}T_{-i}=\dfrac{1}{2}(1-T_{-n-1}-T_{-n+1}).
\end{equation*}%
For a proof of the above formula, see Kuhapatanakul and Sukruan [\ref%
{bib:kuhapatanakul2014}].

Next, we present the formulas which give the summation of the first $n$
Tribonacci and Tribonacci-Lucas sedenions.

\begin{theorem}
The summation formula for Tribonacci and Tribonacci-Lucas sedenions are%
\begin{equation}
\sum\limits_{i=0}^{n}\widehat{T}_{i}=\dfrac{1}{2}(\widehat{T}_{n+2}+\widehat{%
T}_{n}+c_{1})  \label{equation:yaomazxtyhnbvgfc}
\end{equation}%
and%
\begin{equation}
\sum\limits_{i=0}^{n}\widehat{K}_{i}=\dfrac{1}{2}(\widehat{K}_{n+2}+\widehat{%
K}_{n}+c_{2})  \label{equation:mguouresxcdz}
\end{equation}%
respectively, where 
\begin{eqnarray*}
c_{1}
&=&-1-e_{1}-3e_{2}-5e_{3}-9e_{4}-17e_{5}-31e_{6}-57e_{7}-105e_{8}-193e_{9} \\
&&-355e_{10}-653e_{11}-1201e_{12}-2209e_{13}-4063e_{14}-7473e_{15}
\end{eqnarray*}
and 
\begin{eqnarray*}
c_{2}
&=&-6e_{1}-8e_{2}-14e_{3}-28e_{4}-50e_{5}-92e_{6}-170e_{7}-312e_{8}-574e_{9}
\\
&&-1056e_{10}-1842e_{11}-3572e_{12}-6570e_{13}-12084e_{14}-22226e_{15}.
\end{eqnarray*}
\end{theorem}

\textit{Proof.} Using (\ref{equation:bnuyvbcxdszaerm}) and (\ref%
{equation:easzpuyhbnvcdf}), we obtain%
\begin{eqnarray*}
\sum\limits_{i=0}^{n}\widehat{T}_{i}
&=&\sum\limits_{i=0}^{n}T_{i}+e_{1}\sum\limits_{i=0}^{n}T_{i+1}+e_{2}\sum%
\limits_{i=0}^{n}T_{i+2}+...+e_{15}\sum\limits_{i=0}^{n}T_{i+15} \\
&=&(T_{0}+...+T_{n})+e_{1}(T_{1}+...+T_{n+1}) \\
&&+e_{2}(T_{2}+...+T_{n+2})+...+e_{15}(T_{15}+...+T_{n+15}).
\end{eqnarray*}

and 
\begin{eqnarray*}
2\sum\limits_{i=0}^{n}\widehat{T}_{i}
&=&(T_{n+2}+T_{n}-1)+e_{1}(T_{n+3}+T_{n+1}-1-2T_{0}) \\
&&+e_{2}(T_{n+4}+T_{n+3}-1-2(T_{0}+T_{1})) \\
&&\vdots \\
&&+e_{15}(T_{n+17}+T_{n+15}-1-2(T_{0}+T_{1}+...+T_{14})) \\
&=&\widehat{T}_{n+2}+\widehat{T}_{n}+c_{1}
\end{eqnarray*}%
where $%
c_{1}=-1+e_{1}(-1-2T_{0})+e_{2}(-1-2(T_{0}+T_{1}))+...+e_{15}(-1-2(T_{0}+...+T_{14})). 
$ Hence 
\begin{equation*}
\sum\limits_{i=0}^{n}\widehat{T}_{i}=\frac{1}{2}(\widehat{T}_{n+2}+\widehat{T%
}_{n}+c_{1}).
\end{equation*}%
We can compute $c_{1}$ as%
\begin{eqnarray*}
c_{1}
&=&-1-e_{1}-3e_{2}-5e_{3}-9e_{4}-17e_{5}-31e_{6}-57e_{7}-105e_{8}-193e_{9} \\
&&-355e_{10}-653e_{11}-1201e_{12}-2209e_{13}-4063e_{14}-7473e_{15}.
\end{eqnarray*}

This proves (\ref{equation:yaomazxtyhnbvgfc}). Similarly we can obtain (\ref%
{equation:mguouresxcdz}). 
\endproof%

\section{Some Identities For The Tribonacci and Tribonacci-Lucas Sedenions}

In this section, we give identities about Tribonacci and Tribonacci-Lucas
sedenions.

\begin{theorem}
\label{theorem:bnutyvbcdsaerbv}For $n\geq 1,$ the following identities hold:

\begin{description}
\item[(a)] $\widehat{K}_{n}=3\widehat{T}_{n+1}-2\widehat{T}_{n}-\widehat{T}%
_{n-1},$

\item[(b)] $\widehat{T}_{n}+\overline{\widehat{T}_{n}}=2T_{n},$ $\ \widehat{K%
}_{n}+\overline{\widehat{K}_{n}}=2K_{n},$

\item[(c)] $\widehat{T}_{n+1}+\widehat{T}_{n}=\dfrac{\widehat{\alpha }\left(
\alpha +1\right) \alpha ^{n+1}}{(\alpha -\beta )(\alpha -\gamma )}+\dfrac{%
\widehat{\beta }\left( \beta +1\right) \beta ^{n+1}}{(\beta -\alpha )(\beta
-\gamma )}+\dfrac{\widehat{\gamma }\left( \gamma +1\right) \gamma ^{n+1}}{%
(\gamma -\alpha )(\gamma -\beta )},$

\item[(d)] $\widehat{K}_{n+1}+\widehat{K}_{n}=\widehat{\alpha }\left( \alpha
+1\right) \alpha ^{n}+\widehat{\beta }\left( \beta +1\right) \beta ^{n}+%
\widehat{\gamma }\left( \gamma +1\right) \gamma ^{n},$

\item[(e)] $\sum\limits_{i=0}^{n}\binom{n}{i}\widehat{F}_{i}=\dfrac{\widehat{%
\alpha }\alpha (1+\alpha )^{n}}{(\alpha -\beta )(\alpha -\gamma )}+\dfrac{%
\widehat{\beta }\beta (1+\beta )^{n}}{(\beta -\alpha )(\beta -\gamma )}+%
\dfrac{\widehat{\gamma }\gamma (1+\gamma )^{n}}{(\gamma -\alpha )(\gamma
-\beta )},$

\item[(f)] $\sum\limits_{i=0}^{n}\binom{n}{i}\widehat{K}_{i}=\widehat{\alpha 
}(1+\alpha )^{n}+\widehat{\beta }(1+\beta )^{n}+\widehat{\gamma }(1+\gamma
)^{n}.$
\end{description}
\end{theorem}

\textit{Proof. }Since $K_{n}=3T_{n+1}-2T_{n}-T_{n-1}$ (see for example [\ref%
{bib:cerdamorale2018a}]), (a) follows. The others can be established
easily.\ \ 
\endproof%

\begin{theorem}
For $n\geq 0,$ $m\geq 3$, we have

\begin{description}
\item[(a)] $\widehat{T}_{m+n}=T_{m-1}\widehat{T}_{n+2}+(T_{m-2}+T_{m-3})%
\widehat{T}_{n+1}+T_{m-2}\widehat{T}_{n},$

\item[(b)] $\widehat{T}_{m+n}=T_{m+2}\widehat{T}_{n-1}+(T_{m+1}+T_{m})%
\widehat{T}_{n-2}+T_{m+1}\widehat{T}_{n-3},$

\item[(c)] $\widehat{K}_{m+n}=K_{n-1}\widehat{T}_{m+2}+(\widehat{T}_{m+1}+%
\widehat{T}_{m})K_{n-2}+K_{n-3}\widehat{T}_{m+1},$

\item[(d)] $\widehat{K}_{m+n}=K_{m+2}\widehat{T}_{n-1}+(K_{m+1}+K_{m})%
\widehat{T}_{n-2}+K_{m+1}\widehat{T}_{n-3}.$
\end{description}
\end{theorem}

\textit{Proof.} (a) and (d) can be proved by strong induction on $m$ and (c)
can be proved by strong induction on $n.$ For (b), replace $n$ with $n-3$
and $m$ with $m+3$ in (a). 
\endproof%


\begin{thebibliography}{99}
\bibitem{baez2002} \label{bib:baez2002}Baez, J., The octonions, Bull. Amer.
Math. Soc. 39 (2), 145-205, 2002.

\bibitem{biss2008} \label{bib:biss2008}Biss, D.K., Dugger, D., and Isaksen,
D.C., Large annihilators in Cayley-Dickson algebras, Communication in
Algebra, 2008.

\bibitem{bilgici2017} \label{bib:bilgici2017}Bilgici, G., Toke\c{s}er, \"{U}%
,. \"{U}nal, Z., Fibonacci and Lucas Sedenions, Journal of Integer
Sequences, Article 17.1.8, 20, 1-11. 2017.

\bibitem{bruce1984} \label{bib:bruce1984}Bruce, I., A modified Tribonacci
sequence, The Fibonacci Quarterly, 22 : 3, pp. 244--246, 1984.

\bibitem{cariow2013} \label{bib:cariow2013}Cariow, A., Cariowa G., An
Algorithm for Fast Multiplication of Sedenios, Information Proccessing
Letters, Volume 113, Issue, 9, 324-331, 2013.

\bibitem{catalani2002} \label{catalani2002}Catalani, M., Identities for
Tribonacci-related sequences - arXiv preprint,
https://arxiv.org/pdf/math/0209179.pdf math/0209179, 2002.

\bibitem{catarino2018} \label{bib:catarino2018}Catarino, P., k-Pell,
k-Pell--Lucas and modified k-Pell sedenions, Asian-European Journal of
Mathematics, 2018

\bibitem{cawagas2004} \label{bib:cawagas2004}Cawagas, E.R., On the Structure
and Zero Divisors of the Cayley-Dickson Sedenion Algebra, Discussiones
Mathematicae, General Algebra and Applications 24, 251-265, 2004.

\bibitem{cerdamorale2017} \label{bib:cerdamorale2017}Cerda-Morales, G., On a
Generalization for Tribonacci Quaternions, Mediterranean Journal of
Mathematics, 14:239, 1--12, 2017.

\bibitem{cerdamorale2017b} \label{bib:cerdamorale2017b}Cerda-Morales, G.,
The Third Order Jacobsthal Octonions: Some Combinatorial Properties,
arXiv:1710.00602v1, [Math.RA], 2 Oct 2017.

\bibitem{cerdamorale2018a} \label{bib:cerdamorale2018a}Cerda-Morales, G., A
Three-By-Three Matrix Representation of a Generalized Tribonacci sequence,
arXiv:1807.03340v1, [Math.CO], 9 Jul 2018.

\bibitem{cerdamorale2018b} \label{bib:cerdamorale2018b}Cerda-Morales, G.,
The Unifying Formula for All Tribonacci-Type Octonions Sequences and Their
Properties, arXiv:1807.04140v1 [math.RA], [Math.CO], 10 Jul 2018.

\bibitem{choi2013} \label{bib:choi2013}Choi, E., Modular tribonacci Numbers
by Matrix Method, J. Korean Soc. Math. Educ. Ser. B: Pure Appl. Math. Volume
20, Number 3 (August 2013), Pages 207--221, 2013.

\bibitem{cimen2017} \label{bib:cimen2017}\c{C}imen, C., \.{I}pek, A., On
Jacobsthal and Jacobsthal-Lucas Sedenios and Several Identities Involving
These Numbers, Mathematica Aeterna, Vol. 7, No.4, 447-454, 2017.

\bibitem{cimen2017b} \label{bib:cimen2017b}\c{C}imen, C., \.{I}pek, A, On
Jacobsthal and Jacobsthal-Lucas Octonions, Mediterr. J. Math., 14:37, 1-13,
2017.

\bibitem{devbhadra2011} \label{bib:devbhadra2011}Devbhadra, S. V., Some
Tribonacci Identities, Mathematics Today Vol.27(Dec-2011) 1-9, 2011.

\bibitem{elia2001} \label{bib:elia2001}Elia, M., Derived Sequences, The
Tribonacci Recurrence and Cubic Forms, The Fibonacci Quarterly, 39:2, pp.
107-115, 2001.

\bibitem{feinberg1963} \label{bib:feinberg1963}Feinberg, M.,
Fibonacci--Tribonacci, The Fibonacci Quarterly, 1 : 3 (1963) pp. 71--74,
1963.

\bibitem{frontczak2018} \label{bib:frontczak2018}Frontczak, R., Sums of
Tribonacci and Tribonacci-Lucas Numbers, International Journal of
Mathematical Analysis, Vol. 12, No. 1, 19-24, 2018.

\bibitem{gul2018} \label{bib:gul2018}G\"{u}l, K., On k-Fibonacci and k-Lucas
Trigintaduonions, International Journal of Contemporary Mathematical
Sciences, Vol. 13, no. 1, 1 - 10, 2018.

\bibitem{halici2017} \label{bib:halici2017}Halici, S., Karata\c{s}, A., On a
Generalization for Fibonacci Quaternions. Chaos Solitons and Fractals 98,
178--182, 2017.

\bibitem{horadam1963aa} \label{bib:horadam1963aa}Horadam, A. F., Complex
Fibonacci Numbers and Fibonacci quaternions, Amer. Math. Monthly 70,
289--291, 1963.

\bibitem{imaeda2000} \label{bib:imaeda2000}Imaeda, K., Imaeda, M.,
Sedenions: algebra and analysis, Applied Mathematics and Computation, 115,
77-88, 2000.

\bibitem{kecilioglu} \label{bib:kecilioglu}Ke\c{c}ilioglu O, Akku\c{s}, I.,
The Fibonacci Octonions, Adv. Appl. Clifford Algebr. 25, 151--158, 2015.

\bibitem{kilic2008} \label{bib:kilic2008}K\i l\i \c{c}, E., Tribonacci
Sequences with Certain Indices and Their Sums, Ars. Comb., 86, 13-22, 2008.

\bibitem{koplinger2007a} \label{bib:koplinger2007a}K\"{o}plinger, J.,
Signature of gravity in conic sedenions, Applied Mathematics and
Computation, 188, 942-947, 2007.

\bibitem{koplinger2007b} \label{bib:koplinger2007b}K\"{o}plinger, J.,
Gravity and eletromagnetism an conic sedenions, Applied Mathematics and
Computation, 188, 948-953, 2007.

\bibitem{kuhapatanakul2014} \label{bib:kuhapatanakul2014}Kuhapatanakul, K.,
Sukruan, L., The Generalized Tribonacci Numbers With Negative Subscripts,
Integer, 14, 1-6, 2014.

\bibitem{moreno1998} \label{bib:moreno1998}Moreno, G., The zero divisors of
the Cayley-Dickson algebras over the real numbers, Bol. Soc. Mat. Mexicana
(3) 4 , 13-28,1998.

\bibitem{muses1980} \label{bib:muses1980}Muses, C.A., Hypernumber and
quantum field theory with a summary of physically applicable hypernumber
arithmetics and their geometrics, Applied Mathematics and Computation, 6,
63-94, 1980.

\bibitem{okubo1995} \label{bib:okubo1995}Okubo, S., Introduction to
Octonions and Other Non-Associative Algebras in Physics, Cambridge
University Press, Cambridge 1995.

\bibitem{polatli2016} \label{bib:polatli2016}Polatl\i , E., A Generalization
of Fibonacci and Lucas Quaternions, Advances in Applied Clifford Algebras,
26 (2), 719-730, 2016.

\bibitem{lin1988} \label{bib:lin1988}Lin, P. Y., De Moivre-Type Identities
For The Tribonacci Numbers, The Fibonacci Quarterly, 26, pp. 131-134, 1988.

\bibitem{pethe1988} \label{bib:pethe1988}Pethe, S., Some Identities, The
Fibonacci Quarterly, 26, pp. 144--246, 1988.

\bibitem{scott1977} \label{bib:scott1977}Scott, A., Delaney, T., Hoggatt
Jr.,\ V., The Tribonacci sequence, The Fibonacci Quarterly, 15:3, pp.
193--200, 1977.

\bibitem{shannon1977} \label{bib:shannon1977}Shannon, A., Tribonacci numbers
and Pascal's pyramid, The Fibonacci Quarterly, 15:3, pp. 268-275, 1977.

\bibitem{spickerman1981} \label{bib:spickerman1981}Spickerman, W., Binet's
formula for the Tribonacci sequence, The Fibonacci Quarterly, 20,
pp.118--120, 1981.

\bibitem{spickerman1984} \label{bib:spickerman1984}Spickerman, W., Joyner,
R. N., Binets's formula for the Recursive sequence of Order K, The Fibonacci
Quarterly, 22, pp.327-331, 1984.

\bibitem{yalavigi1972} \label{bib:yalavigi1972}Yalavigi, C. C., Properties
of Tribonacci numbers, The Fibonacci Quarterly, 10 : 3, pp. 231--246, 1972.

\bibitem{yilmaz 2014} \label{bib:yilmaz 2014}Yilmaz, N, and N. Taskara, N.,
Tribonacci and Tribonacci-Lucas Numbers via the Determinants of Special
Matrices, Applied Mathematical Sciences, 8, no. 39, 1947-1955, 2014.
\end{thebibliography}
\end{document}